  \magnification 1200

  \catcode `\@=11 \def \@notdefinable {} \catcode `\@=12
  \input miniltx
  \input amssym


  \font \bbfive = bbm5
  \font \bbeight = bbm8
  \font \bbten = bbm10
  \font \rs = rsfs10 \font \rssmall = rsfs10 scaled 833  
  \font \eightbf = cmbx8
  \font \eighti = cmmi8 \skewchar \eighti = '177
  \font \fouri = cmmi5 scaled 800 
  \font \eightit = cmti8
  \font \eightrm = cmr8
  \font \eightsl = cmsl8
  \font \eightsy = cmsy8 \skewchar \eightsy = '60
  \font \eighttt = cmtt8 \hyphenchar \eighttt = -1

  \font \sixi = cmmi6 \skewchar \sixi = '177
  \font \sixrm = cmr6
  \font \sixsy = cmsy6 \skewchar \sixsy = '60
  \font \tensc = cmcsc10

  \scriptfont \bffam = \bbeight
  \scriptscriptfont \bffam = \bbfive
  \textfont \bffam = \bbten

  \newskip \ttglue

  \def \eightpoint {\def \rm {\fam 0 \eightrm }\relax
  \textfont 0 = \eightrm \scriptfont 0 = \sixrm \scriptscriptfont 0 = \fiverm
  \textfont 1 = \eighti \scriptfont 1 = \sixi \scriptscriptfont 1 = \fouri
  \textfont 2 = \eightsy \scriptfont 2 = \sixsy \scriptscriptfont 2 = \fivesy
  \textfont 3 = \tenex \scriptfont 3 = \tenex \scriptscriptfont 3 = \tenex
  \def \it {\fam \itfam \eightit }\relax
  \textfont \itfam = \eightit
  \def \sl {\fam \slfam \eightsl }\relax
  \textfont \slfam = \eightsl
  \def \bf {\fam \bffam \eightbf }\relax
  \textfont \bffam = \bbeight \scriptfont \bffam = \bbfive \scriptscriptfont \bffam = \bbfive
  \def \tt {\fam \ttfam \eighttt }\relax
  \textfont \ttfam = \eighttt
  \tt \ttglue = .5em plus.25em minus.15em
  \normalbaselineskip = 9pt
  \def \MF {{\manual opqr}\-{\manual stuq}}\relax
  \let \sc = \sixrm
  \let \big = \eightbig
  \let \rs = \rssmall
  \setbox \strutbox = \hbox {\vrule height7pt depth2pt width0pt}\relax
  \normalbaselines \rm }

  \def \setfont #1{\font \auxfont =#1 \auxfont }
  \def \withfont #1#2{{\setfont {#1}#2}}


  \def \text #1{{\mathchoice {\hbox {\rm #1}} {\hbox {\rm #1}} {\hbox {\eightrm #1}} {\hbox {\sixrm #1}}}}
  \def \varbox #1{\setbox 0\hbox {$#1$}\setbox 1\hbox {$I$}{\ifdim \ht 0< \ht 1 \scriptstyle #1 \else \scriptscriptstyle #1 \fi }}

  \def \rsbox #1{{\mathchoice {\hbox {\rs #1}} {\hbox {\rs #1}} {\hbox {\rssmall #1}} {\hbox {\rssmall #1}}}}
  \def \mathscr #1{\rsbox {#1}}

  \font \rstwo = rsfs10 scaled 900
  \font \rssmalltwo = rsfs10 scaled 803


  \def \TRUE {Y}
  \def \FALSE {N}

  \def \ifundef #1{\expandafter \ifx \csname #1\endcsname \relax }

  \def \undefrule {\kern 2pt \vrule width 2pt height 5pt depth 0pt \kern 2pt}
  \def \UndefLabels {}
  \def \possundef #1{\ifundef {#1}\undefrule {\eighttt #1}\undefrule
    \global \edef \UndefLabels {\UndefLabels #1\par }
  \else \csname #1\endcsname \fi }


  \newcount \secno \secno = 0
  \newcount \stno \stno = 0
  \newcount \eqcntr \eqcntr = 0

  \ifundef {showlabel} \global \def \showlabel {\FALSE } \fi
  \ifundef {auxwrite} \global \def \auxwrite {\TRUE } \fi
  \ifundef {auxread} \global \def \auxread {\TRUE } \fi

  \def \newcommand #1#2{\global \edef #1{#2}} 

  \def \define #1#2{\global \expandafter \edef \csname #1\endcsname {#2}}
  \long \def \error #1{\medskip \noindent {\bf ******* #1}}
  \def \fatal #1{\error {#1\par Exiting...}\end }

  \def \advseqnumbering {\global \advance \stno by 1 \global \eqcntr =0}

  \def \current {\ifnum \secno = 0 \number \stno \else \number \secno \ifnum \stno = 0 \else .\number \stno \fi \fi }


  \def \rem #1{\vadjust {\vbox to 0pt{\vss \hfill \raise 3.5pt \hbox to 0pt{ #1\hss }}}}
  \font \tiny = cmr6 scaled 800
  \def \deflabel #1#2{\relax
    \if \TRUE \showlabel \rem {\tiny #1}\fi
    \ifundef {#1PrimarilyDefined}\relax
      \define {#1}{#2}\relax
      \define {#1PrimarilyDefined}{#2}\relax
      \if \TRUE \auxwrite \immediate \write 1 {\string \newlabe l {#1}{#2}}\fi
    \else
      \edef \old {\csname #1\endcsname }\relax
      \edef \new {#2}\relax
      \ifx \old \new \else \fatal {Duplicate definition for label ``{\tt #1}'', already defined as ``{\tt \old }''.}\fi
      \fi }

  \def \label #1 {\deflabel {#1}{\current }}

  \def \lbldeq #1 $$#2$${\if \InsideProof N \advseqnumbering
    \edef \lbl {\current }\else
    \global \advance \eqcntr by 1
    \edef \lbl {\current .\number \eqcntr }\fi
    \deflabel {#1}{\lbl }
    $$
    #2
    \eqno {(\lbl )}
    $$}

  \def \split #1.#2.#3.#4;{\global \def \parone {#1}\global \def \partwo {#2}\global \def \parthree {#3}\global \def
\parfour {#4}}
  \def \NA {NA}
  \def \ref #1{\split #1.NA.NA.NA;(\possundef {\parone }\ifx \partwo \NA \else .\partwo \fi )}


  \newcount \bibno \bibno = 0

  \def \Bibitem #1 #2; #3; #4 \par {\smallbreak
    \global \advance \bibno by 1
    \item {[\possundef {#1}]} #2, {``#3''}, #4.\par
    \ifundef {#1PrimarilyDefined}\else
      \fatal {Duplicate definition for bibliography item ``{\tt #1}'',
      already defined in ``{\tt [\csname #1\endcsname ]}''.}
      \fi \ifundef {#1}\else \edef \prevNum {\csname #1\endcsname } \ifnum \bibno =\prevNum \else \error {Mismatch
        bibliography item ``{\tt #1}'', defined earlier (in aux file ?) as ``{\tt \prevNum }'' but should be ``{\tt
        \number \bibno }''.  Running again should fix this.}  \fi \fi
    \define {#1PrimarilyDefined}{#2}\relax
    \if \TRUE \auxwrite \immediate \write 1 {\string \newbi b {#1}{\number \bibno }}\fi }

  \def \jrn #1, #2 (#3), #4-#5;{{\sl #1}, {\bf #2} (#3), #4--#5}
  \def \Article #1 #2; #3; #4 \par {\Bibitem #1 #2; #3; \jrn #4; \par }

  \def \references {\begingroup \bigbreak \eightpoint \centerline {\tensc References} \nobreak \medskip \frenchspacing }


  \catcode `\@ =11
  \def \citetrk #1{{\bf \possundef {#1}}} 
  \def \c@ite #1{{\rm [\citetrk {#1}]}}
  \def \sc@ite [#1]#2{{\rm [\citetrk {#2}\hskip 0.7pt:\hskip 2pt #1]}}
  \def \du@lcite {\if \pe@k [\expandafter \sc@ite \else \expandafter \c@ite \fi }
  \def \cite {\futurelet \pe@k \du@lcite }
  \catcode `\@ =12


  \def \Headlines #1#2{\nopagenumbers
    \headline {\ifnum \pageno = 1 \hfil
    \else \ifodd \pageno \tensc \hfil \lcase {#1} \hfil \folio
    \else \tensc \folio \hfil \lcase {#2} \hfil
    \fi \fi }}

  \def \title #1{\medskip \centerline {\withfont {cmbx12}{\ucase {#1}}}}

  \def \Subjclass #1#2{\footnote {\null }{\eightrm #1 \eightsl Mathematics Subject Classification: \eightrm #2.}}

  \long \def \Quote #1\endQuote {\begingroup \leftskip 35pt \rightskip 35pt \parindent 17pt \eightpoint #1\par \endgroup }
  \long \def \Abstract #1\endAbstract {\vskip 1cm \Quote \noindent #1\endQuote }

  \def \Note #1{\footnote {}{\eightpoint #1}}
  \def \Date #1 {\Note {\it Date: #1.}}

  \newcount \auxone \newcount \auxtwo \newcount \auxthree
  \def \currenttime {\auxone =\time \auxtwo =\time \divide \auxone by 60 \auxthree =\auxone \multiply \auxthree by 60
\advance
    \auxtwo by -\auxthree \ifnum \auxone <10 0\fi \number \auxone :\ifnum \auxtwo <10 0\fi \number \auxtwo }
  \def \today {\ifcase \month \or January\or February\or March\or April\or May\or June\or July\or August\or September\or
    October\or November\or December\fi { }\number \day , \number \year }
  
  \def \hojeExtenso {\number \day \ de \ifcase \month \or janeiro\or fevereiro\or mar\c co\or abril\or maio\or junho\or
julho\or
    agosto\or setembro\or outubro\or novembro\or decembro\fi \ de \number \year }

  \def \part #1#2{\vfill \eject \null \vskip 0.3\vsize
    \withfont {cmbx10 scaled 1440}{\centerline {PART #1} \vskip 1.5cm \centerline {#2}} \vfill \eject }


  \def \Case #1:{\medskip \noindent {\tensc Case #1:}}

  \def \fix {\smallskip \noindent $\blacktriangleright $\kern 12pt}

  \def \ucase #1{\edef \auxvar {\uppercase {#1}}\auxvar }
  \def \lcase #1{\edef \auxvar {\lowercase {#1}}\auxvar }

  \def \emph #1{{\it #1}\/}

  \def \state #1 #2\par {\begingroup \def \InsideProof {Y} \medbreak \noindent \advseqnumbering {\bf \current .\enspace
#1.\enspace \sl #2\par }\medbreak \endgroup }

  \def \lstate #1 #2 #3\par {\state {#2} \label {#1} #3\par}

  \def \definition #1\par {\state Definition \rm #1\par }




  \def \Proof {\if \InsideProof N \else \fatal {Opening Proof block before closing the last one} \fi
    \global \def \InsideProof {Y}
    \medbreak \noindent {\it Proof.\enspace }}

  \def \endProof {\if \InsideProof Y \hfill $\endproofmarker $ \looseness = -1 \fi
    \medbreak
    \if \InsideProof N \fatal {Closing Proof block before opening it} \fi
    \global \def \InsideProof {N}}

  \def \InsideProof {N}


  \def \explica #1#2{\mathrel {\buildrel \hbox {\sixrm #1} \over #2}}
  \def \explain #1#2{\explica {\ref {#1}}{#2}} 
  
  \def \=#1{\explain {#1}{=}}

  \newcount \fnctr \fnctr = 0
  \def \fn #1{\global \advance \fnctr by 1
    \edef \footnumb {$^{\number \fnctr }$}\relax
    \footnote {\footnumb }{\eightpoint #1\par \vskip -10pt}}

  \def \text #1{\mathchoice {\hbox {#1}} {\hbox {#1}} {\hbox {\eightrm #1}} {\hbox {\sixrm #1}}}


  \def \item #1{\par \noindent \kern 1.1truecm\hangindent 1.1truecm \llap {#1\enspace }\ignorespaces }
  
  \def \Item #1{\smallskip \item {{\rm #1}}}

  \newcount \zitemno \zitemno = 0
  \def \izitem {\global \zitemno = 0}

  \def \zitemplus {\global \advance \zitemno by 1 \relax }
  \def \rzitem {\romannumeral \zitemno }
  \def \rzitemplus {\zitemplus \rzitem }
  \def \zitem {\Item {{\rm (\rzitemplus )}}}

  \def \lzitem #1 {\Item {{\rm (\rzitemplus )}} {\deflabel {#1}{\current .\rzitem }{\def \showlabel {\FALSE }\deflabel
{Local#1}{\rzitem }}}}

  \newcount \nitemno \nitemno = 0
  
  \def \nitem {\global \advance \nitemno by 1 \Item {{\rm (\number \nitemno )}}}

  \newcount \aitemno \aitemno = -1
  \def \boxlet #1{\hbox to 6.5pt{\hfill #1\hfill }}
  \def \iaitem {\aitemno = -1}
  \def \aitemconv {\ifcase \aitemno a\or b\or c\or d\or e\or f\or g\or h\or i\or j\or k\or l\or m\or n\or o\or p\or q\or
    r\or s\or t\or u\or v\or w\or x\or y\or z\else zzz\fi }
  \def \aitem {\global \advance \aitemno by 1\Item {(\boxlet \aitemconv )}}


  \def \deflabeloc #1#2{\deflabel {#1}{\current .#2}{\edef \showlabel {\FALSE }\deflabel {Local#1}{#2}}}
  \def \lbldzitem #1 {\zitem \deflabeloc {#1}{\rzitem }}
  \def \lbldaitem #1 {\aitem \deflabeloc {#1}{\aitemconv }}
  \def \aitemmark #1 {\deflabel {#1}{\aitemconv }}

  \def \zitemmark #1 {\deflabel {#1}{\current .\rzitem }{\def \showlabel {\FALSE }\deflabel {Local#1}{\rzitem }}}


  \def \mathbb #1{{\bf #1}}
  \def \frac #1#2{{#1\over #2}}

  \def \<{\left \langle \vrule width 0pt depth 0pt height 8pt }
  \def \>{\right \rangle }
  \def \({\big (}
  \def \){\big )}
  
  \def \and {\mathchoice {\hbox {\quad and \quad }} {\hbox { and }} {\hbox { and }} {\hbox { and }}}

  \def \IFF {\kern 7pt\Leftrightarrow \kern 7pt}
  \def \IMPLY {\kern 7pt \Rightarrow \kern 7pt}
  \def \for #1{\mathchoice {\quad \forall \,#1} {\hbox { for all } #1} {\forall #1}{\forall #1}}
  \def \endproofmarker {\square }
  \def \"#1{{\it #1}\/} 
  
  \def \*{\otimes }
  \def \caldef #1{\global \expandafter \edef \csname #1\endcsname {{\cal #1}}}
  \def \mathcal #1{{\cal #1}}
  \def \bfdef #1{\global \expandafter \edef \csname #1\endcsname {{\bf #1}}}
  \bfdef N \bfdef Z \bfdef C \bfdef R
  \def \exists {\mathchar "0239\kern 1pt }
  \def \labelarrow #1{\setbox 0\hbox {\ \ $#1$\ \ }\ {\buildrel \textstyle #1 \over {\hbox to \wd 0 {\rightarrowfill }}}\ }
  \def \subProof #1{\medskip \noindent #1\enspace }
  \def \itmProof (#1) {\subProof {(#1)}}
  \def \itemImply #1#2{\subProof {#1$\Rightarrow $#2}}
  \def \itmImply (#1) > (#2) {\itemImply {(#1)}{(#2)}}


  \if \TRUE \auxread
    %
    \IfFileExists {\jobname .aux}{\input \jobname .aux}{\null } \fi
  \if \TRUE \auxwrite \immediate \openout 1 \jobname .aux \fi


  \def \close {\ifx \empty \UndefLabels \else
    \message {*** There were undefined labels ***} \medskip \noindent
    ****************** \ Undefined Labels: \tt \par \UndefLabels \fi
    \if \TRUE \auxwrite \closeout 1 \fi
    \par \vfill \supereject \end }

  %
  %

  \def \startsection #1 \par
    {\goodbreak \bigbreak
    \begingroup
    \global \edef \secname {#1}\relax
    \global \advance \secno by 1
    \stno = 0}

  \def \sectiontitle \par
    {\noindent {\bf \number \secno .\enspace \secname .}
    \nobreak \medskip
    \noindent }

  \def \endsection {\endgroup }

\def \G{{\cal G}}
\def \rsboxtwo #1{{\mathchoice {\hbox {\rstwo #1}} {\hbox {\rstwo #1}} {\hbox {\rssmalltwo #1}} {\hbox {\rssmalltwo #1}}}}
\def \Lbsys {\rsboxtwo {L}}
\def \Lb {\Lbsys \kern 1pt}
\def \Lbpar {\Lbsys \kern 2pt}
\def \P{\rsboxtwo{P}\, }
\def \Ep{\rsboxtwo{EP}\, (B)}

\def \red {_{\hbox {\sixrm red}}}
\def \NAB{N(A, B)}
\def \Gz{\G^{(0)}}
\def \eg{\'etale groupoid}
\def \teg{twisted \eg}

\def \tpeg{twisted, principal, \eg}
\def\Nuc{{L}}
\long \def \Abstract #1\endAbstract {\Quote \noindent #1\endQuote }

\def\EP#1#2#3{\cite[#1 #3]{ExelPitts}}


\def \LocProj{3.6}
\def \DefBlackIdeal{11.4}


\def \CorrespondenceExtendedStates{3.20}
\def \BetaDomain{4.4.iii}
\def \BetaIdentity{4.7}
\def \Beta{4.4.ii}
\def \CondexForState{7.3}
\def \NorGenCx{7.5}
\def \DefFreePts{9.1}
\def \DefineCano{20.1}
\def \DefineRelTriv{18.1}
\def \DescribeGray{15.3}
\def \EbContinuous{8.4}
\def \IntroBlack{11.2}
\def \IsotropyOneDim{9.3}
\def \MainThree{20.6}
\def \MainThree{20.6}
\def \MainThree{20.6}
\def \ThreeValues{9.4}
\def \TwoAlgebras{14.8}
\def \IdIsoAlgsFull{14.12.ii}
\def \CharacFreePts{14.14}


\def\.{ \kern 1pt}
\font \tit = cmbx10
\font\itt = cmtt10 scaled 700
\font\iss = cmss10 scaled 800

\def\.{ \kern 2pt}
\centerline {\tit ON\.KUMJIAN'S\.C*-DIAGONALS\.AND}
\smallskip
\centerline {\tit THE\.OPAQUE\.IDEAL}

\bigskip
\centerline {\tensc
  R. Exel\footnote {$^{\ast }$}{\eightrm Universidade Federal de Santa Catarina. Partially supported by CNPq.}
  }

\footnote {}{\eightrm October 16,  2021.}

\Subjclass {2010}{46L55, 46L30}

\Abstract
  We characterize exotic C*-algebras of {\tpeg}s, together with the abelian subalgebra associated to the unit space, as
precisely being the inclusions ``$A\subseteq B$'' of C*-algebras in which $A$ is abelian, regular, and satisfies the extension
property (pure states extend uniquely to $B$).
When $B$ is moreover nuclear, we deduce that the corresponding opaque ideal is trivial.  As an
application, we give a streamlined characterization of Kumjian's C*-diagonals as the regular abelian subalgebras
satisfying the extension property with vanishing opaque ideal.
 \endAbstract

\startsection Introduction

\sectiontitle

Inspired by Feldman and Moore's fundamental paper \cite{FeldmanMoore} on Cartan subalgebras of von Neumann algebras,
Kumjian \cite{Kumjian} introduced the notion of \emph{diagonals} in the context of C*-algebras and proved that every
inclusion
  \lbldeq DaInclusion
  $$
  A\subseteq B
  $$
  of an
  abelian C*-algebra $A$ in another C*-algebra $B$, satisfying suitable hypotheses, is necessarily modeled by a {\tpeg}.

The first main goal of the present paper is to prove a souped up version of Kumjian's result, based on my recent work
\cite{ExelPitts} with D.~Pitts.
  The plan is to strip the hypotheses of \cite[Theorem 3.1]{Kumjian} to a bare minimum, while retaining its conclusion,
except that the modeling will be done with a possibly exotic groupoid C*-algebra, rather than the reduced version
adopted in \cite{Kumjian}.

Besides requiring the indispensable condition that \ref{DaInclusion} be a regular inclusion (Definition
\ref{DefReg.iii} below), Kumjian assumes the existence of a faithful conditional expectation whose kernel is
required to be spanned by
the so called \emph{free} normalizers, namely normalizers squaring to zero.
In an important intermediate result \cite[Proposition 1.4]{Kumjian}  Kumjian shows that C*-diagonals satisfy the \emph{extension property}, meaning that each pure state
of $A$ admits a unique extension to a state on $B$ (see Definition \ref{DefEP} below for more details).

We in turn adopt the extension property as the sole condition (besides regularity) in our main result, namely Theorem
\ref{Main}, where we prove that such inclusions are precisely the ones arising from exotic C*-algebras of {\tpeg}s.

The proof of our main result is achieved as a direct application of \EP{Theorem}{20.6}{\MainThree}, a result that
characterizes exotic
groupoid C*-algebras in full generality.  In order to make use of it, all we need is to provide a
\emph{canonical state} \EP{Definition}{20.1}{\DefineCano} relative to every point of the spectrum of $A$.  This is in turn obtained
as a consequence of two well known results, namely that inclusions satisfying the extension property are maximal abelian
\cite[page 385]{KadisonSinger}, and automatically posses a conditional expectation \cite[Theorem 2.2]{Archbold} (albeit not
necessarily a faithful one).  Since the proofs of these two results are very short, and often only done in the unital case, we
include them for the convenience of the reader.

Our second main  goal, actually the motivation for writing the present paper, is the study of the opaque ideal
for regular inclusions.  We first show that the regular inclusions
with a trivial opaque ideal, and  satisfying the extension property, are precisely
Kumjian's C*-diagonals.  This result generalizes \cite[Theorem 4.8]{DonsigPitts} to the non-unital case.
  Furthermore, assuming the extension property, we give an affirmative answer to a question raised in \cite[Question
4.2]{EPZ}, by proving that the opaque ideal is trivial when the containing algebra is nuclear.

Last but not least, I'd like to thank
  mathoverflow
  user {\iss Darth Vader} for bringing
up Archbold's paper in a comment to a question of mine \cite{DarthVader}, which triggered a long chain of ideas, eventually
culminating with the present paper.  Thanks are also due to David Pitts for several helpful suggestions.
\endsection

\startsection Basic facts

\sectiontitle

\fix Throughout this section we shall fix a C*-algebra $B$, and we shall let
  $$
  A\subseteq B
  $$
  be a fixed closed *-subalgebra.  In our main
result, below, we will assume that $A$ is abelian, so the reader is welcome to assume this from now on, although in a few
places commutativity of $A$ is not strictly required.  When $A$ is explicitly assumed to be abelian, the spectrum of
$A$ will be denoted by $X$, so that $A$ is *-isomorphic to $C_0(X)$ by Gelfand's Theorem.

\definition \label DefEP
One says that
$A$ has the \emph{extension property}
relative to $B$, provided:
  \izitem
  \zitem every pure state of $A$ admits a unique extension to
a state on $B$, and
  \zitem $A$ contains an approximate identity for $B$.

Many variations of the above notion appear in the literature.  Quite often
  [\citetrk {KadisonSinger},\citetrk {JAnderson},\citetrk{ArchboldBunceGregson}]
  it is required that $B$ be unital and that $A$
contain the unit of $B$, which obviously implies point (ii) above.  In \cite{Archbold}, on the other hand, point (ii) is
replaced by the equivalent requirement that no pure state of $B$ annihilates $A$ (the equivalence follows from the
Krein-Milman Theorem and \cite[Lemma 2.32]{Akemann}).

A nice conceptual equivalent formulation is as follows:  denoting by
$\P(A)$ the set of all pure states on $A$, and by $\P_0(A)=\P(A)\cup \{0\}$, and similarly for $B$, then the extension property may
be characterized by saying that every member of $\P_0(A)$ admits a unique extension to a member of $\P_0(B)$.  The catch
is that the zero linear functional on $A$ admits a unique extension to a member of $\P_0(B)$ iff no pure state of $B$ restricts
to zero on $A$.

The following notion was introduced by Kumjian in \cite{Kumjian}:

\definition \label DefReg
  \izitem
  \zitem An element $b\in B$   is said to be a \emph{normalizer} for $A$, if
  $$
  b^*ab\subseteq  A, \and   bab^*\subseteq  A.
  $$
  \zitem The set of all normalizers is denoted by $\NAB$.
  \zitem One says that $A$ is a regular subalgebra of $B$ when  \ref{DefEP.ii} holds, and $\NAB$ spans a dense subspace of $B$.

\state Remarks \rm
  \iaitem
  \aitem Observe that,  assuming \ref{DefEP.ii},  one has that  $n^*n\in A$  for every $n$ in $\NAB$.  This is because if
$\{u_i\}_i$ is an approximate unit for $B$ contained in $A$,  then
  $$
  n^*n = \lim _{i\to \infty } n^*u_in \in  A.
  $$
  \aitem
  On the other hand, if we assume that $\NAB$ spans $B$, and that $n^*n$ lie in $A$ for every normalizer $n$, then
\ref{DefEP.ii} follows.  Indeed, taking any normalizer $n$, we have
  $$
  n = \lim_{k\to \infty }n(n^*n)^{1/k} \in  \overline{BA},
  $$
  from where it follows that $B=\overline{BA}$, and from this it is easily seen that any approximate unit for $A$ is also an
approximate unit for $B$.
  \aitem
This said, it would perhaps be slightly more elegant to strengthen the definition of normalizers by adding the requirement that
$n^*n$, and perhaps also $nn^*$, lie in $A$.    Should this be done, the above definition of a regular subalgebra could be
streamlined by requiring
only that $\NAB$ span $B$, as \ref{DefEP.ii} would be automatic.
  \aitem
We should also remark that D.~Pitts  has recently shown  \cite{Pitts} that, when $A$ is a maximal abelian subalgebra of $B$, then
the fact that $\NAB$ span $B$ (with the standard definition of normalizers) alone implies \ref{DefEP.ii}.  The reader
should however be warned that, unlike \ref{DefReg.iii},
 Pitts defines regularity without condition \ref{DefEP.ii}.


In order to prove our main result we need two auxiliary facts.  The first one refers to the existence of a conditional
expectation from $B$ to $A$,  a fact that has been proved by
  Anderson \cite[Theorem 3.4]{JAnderson} for $A$ unital and abelian, and generalized by
  Archbold \cite[Theorem 2.2]{Archbold} for some situations where $A$ does not need to be commutative.
  In case $A$ is abelian and regular the proof is much shorter, so it is perhaps worth spelling it out.

\state Lemma \label LemmaCondexp
  {\rm [\citetrk {JAnderson},\citetrk{Archbold}]}
  Assume that $A$ is abelian,  regular,  and satisfies the extension
property.  For every $x$ in $X$ (the spectrum  of $A$),
let
  $\varphi _x$ be the pure state defined
  by
  $$
  \varphi _x(a) = a(x), \for a\in  A,
  $$
  and let $\psi _x$ be the unique state on  $B$ extending $\varphi _x$.
  Moreover, for each $b$ in $B$, consider the scalar-valued function $E(b)$ defined  on $X$ by
  $$
  E(b)|_x = \psi _x(b), \for x\in  X.
  $$
  Then
  \izitem
  \zitem $E(b)$ is continuous on $X$,
  \zitem $E(b)$ belongs to $A$ (that is, in case $X$ is not compact,  $E(b)$  vanishes at $\infty $),
  \zitem The ensuing map $E:B\to A$ is the unique conditional expectation from $B$ to $A$.

\Proof
  The continuity of $E(b)$ follows immediately from \EP{Proposition}{8.4}{\EbContinuous}, once we realize that the set  $F_b$ mentioned
there coincides with $X$, and that $\varepsilon _b=E(b)$.

It is also clear that $E(b)$ is bounded with
$\|E(b)\|_\infty \leq \|b\|$, hence $E$ defines a contractive,  positive map
  $$
  E:B\to C^b(X),
  $$
  where $C^b(X)$ stands for the C*-algebra of all bounded, continuous, complex functions on $X$.
If $b\in  B$ and $a\in  A$, then
  $$
  E(ab) = aE(b) = E(ba),
  $$
by \EP{Lemma}{7.3}{\CondexForState}.  Therefore,
letting $\{u_i\}_i$ be an approximate unit for $B$ contained in $A$, we have that
  $$
  E(b) = \lim_{i\to \infty } E(u_ib) = \lim_{i\to \infty } u_iE(b) \in  C_0(X) =A,
  $$
  proving (ii).

Regarding (iii), it is now clear that $E$ is a conditional expectation from $B$ to $A$.
To show that $E$ is unique,  suppose that $F$ is another one.  Then, for
every $x$ in $X$, we
may define  a state $\rho _x$ on $B$ by
  $$
  \rho _x:b\in  B \mapsto  F(b)|_x\in {\bf C}.
  $$
  Evidently $\rho _x$ extends $\varphi _x$, so $\rho _x=\psi _x$ by uniqueness.  Therefore, for every $b$ in $B$, and every $x$ in $X$, we have
that
  $$
  F(b)|_x = \rho _x(b) = \psi _x(b) =   E(b)|_x,
  $$
  proving that $F=E$.
  \endProof

The second auxiliary tool needed to prove our main result is the fact that the extension property implies maximal
commutativity.  In the case of unital algebras, this fact is implicitly mentioned in the introduction to \cite{KadisonSinger},
  incidentally the same paper where the famous Kadison-Singer problem was first posed.
  This  is also stated in
\cite[Corollary 2.7]{ArchboldBunceGregson},  where \cite{KadisonSinger} is referenced for a proof.  In what follows we give a proof of this result for the sake of
completeness, because the available references apparently only deal with the unital case, and also because the proof is
really very short.

\state Lemma \label MaxCommut \cite{KadisonSinger}
  If $A$ is abelian and satisfies the extension property relative to $B$, then $A$ is
maximal abelian.

\Proof
If $b\in B$ commutes with $A$, we must prove that $b\in  A$.  By decomposing $b$ in its real and imaginary parts, we may
assume that $b$ is self-adjoint.
The closed *-algebra $C$ of $B$  generated by $A\cup \{b\}$ is therefore abelian, so we may
write $C=C_0(\Omega )$, for some locally compact space $\Omega $, while the inclusion
  $$
  C_0(X)=A\ \hookrightarrow \ C=C_0(\Omega )
  $$
  is necessarily implemented by a proper, surjective  map
  $h:\Omega \to X$
  (notice that \ref{DefEP.ii} is crucial for such a map to exist).  Observe  that  $h$ must be one-to-one because, if
$s$ and $t$ are distinct points in $\Omega $ such that $h(s)=h(t)=:x$, then
  both $\varphi _s$ and $\varphi _t$ extend $\varphi _x$ from $A$ to $C$.  Further state extensions of    $\varphi _s$ and $\varphi _t$ to $B$ would
then violate the extension property.  This proves that $h$ is bijective. Since $\Omega $ is locally compact and $h$ is
proper, it follows that $h$ is a homeomorphism, so $A=C$, and hence $b\in  A$.
\endProof

\endsection

\startsection Exotic C*-algebras of {\tpeg}s

\def\norm{\|\kern1pt{\cdot}\kern1pt\|}
\def \max {_{\hbox {\sixrm max}}}

\sectiontitle

Recall that,  given any {\teg} $(\G, \Lbpar)$, the \emph{full} and \emph{reduced} C*-algebras of $(\G,
\Lbpar)$,  denoted
  $$
  C^*(\G,\Lbpar)  \and  C^*\red(\G,\Lbpar),
  $$
  are respectively defined to be the completions of $C_c(\G,\Lbpar)$ under the \emph{maximal} and \emph{reduced} norms, namely
  $$
  \norm\max \and \norm\red.
  $$
  If we  instead complete $C_c(\G,\Lbpar)$ under some C*-norm $\mu $ satisfying
  \lbldeq ExoticNorm
  $$
  \|f\|\max \geq  \mu (f) \geq   \|f\|\red, \for f \in  C_c(\G,\Lbpar),
  $$
  the resulting object is often denoted $C^*_\mu (\G,\Lbpar)$, and
called an
  \emph{exotic}\fn{Some authors use of the word \emph{exotic} only in case  $\mu $ is distinct
from either the full or reduced norms.  Nevertheless,  as should be clear from the above,  we shall not make that distinction.}
  groupoid C*-algebra.

Our first main result is a complete characterization of exotic C*-algebras of {\tpeg}s.

\state Theorem \label Main
  Given a C*-algebra $B$ and a closed *-subalgebra $A\subseteq B$, the following are equivalent:
  \iaitem
  \aitem $A$ is abelian, regular, and satisfies the extension property,
  \aitem $B$ is isomorphic to an exotic C*-algebra for a {\tpeg} $(\G,\Lbpar)$,
  via an isomorphism carrying $A$ onto $C_0\big (\Gz\big )$.

\Proof
  Assuming (a)
  we will prove (b)  by means of a direct application of  \EP{Theorem}{20.6}{\MainThree}.  For this we let
  $N = \NAB $,
  and we will prove that the state  $\psi _x$ (introduced in the proof of \ref{LemmaCondexp})
is an \emph{$N$-canonical state} \EP{Definition}{20.1}{\DefineCano} relative to $x$, for every $x$ in $X$.
In  other words,   we must show that $\psi _x$ vanishes on any element $n\in N$, except possibly at those for which
$x$ is \emph{trivial   relative to $n$} \EP{Definition}{18.1}{\DefineRelTriv}, meaning that there exists $v$ in $A$ such that $v(x)=1$,
and $nv\in A$.

Seen from the contra-positive point of view,  our task is to show
that, for every $x$ in $X$, and every $n$ in $N$, one has that
  $$
  \psi _x(n) \neq 0 \IMPLY x \text { is trivial relative to } n.
  $$

  Given that $\psi _x(n) \neq 0$, we deduce from \ref{LemmaCondexp.i} that there exists a neighborhood $V$ of $x$ such that
$\psi _y(n)\neq 0$, for every $y$ in $V$.
  The extension property precisely says that every point of $X$ is free in the sense of
  \EP{Definition}{9.1}{\DefFreePts},
  so the isotropy algebra $B(y)$ is isomorphic to ${\bf C}$ by
  \EP{Proposition}{9.3}{\IsotropyOneDim}.
  Moreover, identifying $B(y)$ with ${\bf C}$, we have that
  $$
  E_y(n) = \psi _y(n),
  $$
  by \EP{Proposition}{9.4}{\ThreeValues},
  where
  $E_y:B\to B(y)$
  is the corresponding localizing projection
  \EP{Definition}{3.6}{\LocProj}
  (for the pair of ideals $(J_y, J_y)$).
  From
  $
  \psi _y(n)  \neq 0,
  $
  and employing
  \EP{Proposition}{7.5.iii}{\NorGenCx.iii},
  we deduce
  that $n$ lies in $N_{y, y}$, that is, $y$ belongs to the domain of $\beta _n$ (see \EP{Proposition}{4.4.ii}{\Beta}), and
$\beta _n(y)=y$.  In other words, $\beta _n$ coincides with the identity map on $V$.

Choose any $v$ in $A$ such that $v(x)=1$, and $v=0$ outside $V$.  Setting $m:= nv$, it follows from
\EP{Proposition}{4.4.iii}{\BetaDomain} that $x$ lies in the domain of $\beta _m$, and that $\beta _m$ is the identity on its
domain.  We may then invoke \EP{Proposition}{4.7}{\BetaIdentity} to conclude that $m$ lies in the commutant $A'$, hence
also in $A$ by \ref{MaxCommut}.  This shows that $nv\in A$, so we see that $x$ is trivial relative to $n$, as desired.  An
application of \EP{Theorem}{20.6}{\MainThree} then takes care of everything, except for the fact that $\G$ is principal.

To see that $\G$ is principal,
  let $B_1=C^*\red(\G,\Lbpar)$.
  Since $\mu $ dominates the reduced norm, $B_1$ is a quotient of $B$, and as $A$ satisfies the extension property relative
to $B$, the same is evidently true with respect to $B_1$.  By \EP{Corollary}{14.14}{\CharacFreePts} (adopting the second
alternative in \EP{}{14.8}{\TwoAlgebras}), we then deduce that every isotropy group $\G(x)$ is trivial, so $\G$ is
principal.

\bigskip
  Conversely,  assuming (b), it is immediate that $A$ is abelian and regular, so it suffices to prove the extension property.
Considering the natural epimorphism
  $$
  \pi : C^*(\G,\Lbpar) \to   C^*_\mu (\G,\Lbpar),
  $$
  arising from the first inequality in \ref{ExoticNorm}, observe that $\pi $ restricts to the identity map on the
respective copies of $C_0\big (\Gz\big )$.  Therefore, for each $x\in X= \Gz$, there are at least as many extensions of
$\varphi _x$ to $C^*(\G,\Lbpar)$, as there are to $C^*_\mu (\G,\Lbpar)$.  This said, it is enough to prove the extension property
under the assumption that $B=C^*(\G,\Lbpar)$.

Given any $x$ in $X$, we have by
\EP{Theorem}{14.12.i}{\IdIsoAlgsFull}
  that the isotropy algebra $B(x)$ is isomorphic to the (twisted) group C*-algebra of the isotropy group $\G(x)$.
However this group is trivial due to the fact that
$\G$ is principal, whence $B(x)={\bf C}$.

Using   \EP{Proposition}{9.3}{\IsotropyOneDim} we then conclude that every $x$ in $X$ is free, meaning that $\varphi _x$ admits a
unique extension to a state on $B$, according to \EP{Definition}{9.1}{\DefFreePts}.
\endProof

We should point out that one could attempt to prove \ref{Main} using a different strategy, based on a result by Donsig
and Pitts \cite[Theorem 4.8]{DonsigPitts}, where it is shown that the left kernel of the conditional expectation $E$
given by \ref{LemmaCondexp}, namely
  $$
  \Nuc =\{b\in  B: E(b^*b)=0\},
  $$
  is a two-sided ideal in $B$, trivially intersecting $A$.
The quotient $B/\Nuc $ therefore contains a faithful, regular
copy of $A$.  Donsig and Pitts go on to prove that $A$ is then a
C*-diagonal in $B/\Nuc $, so Kumjian's Theorem may be applied to
model $B/\Nuc $ as the reduced C*-algebra of a {\tpeg} $(\G, \Lbpar)$.

One would then need to show that the C*-algebra of $(\G, \Lbpar)$, namely
the twisted Weyl groupoid for $B/\Nuc $, admits a representation in $B$ (as opposed to $B/\Nuc $), which would
then lead to a sequence of *-homomorphisms
  $$
  C^*(\G,\Lbpar) \to  B \to  B/\Nuc  \simeq   C^*\red(\G,\Lbpar),
  $$
  from where it would follow that $B$ is an exotic version of   $C^*(\G,\Lbpar)$, as desired.  We don't presently have
a clear strategy for providing such a representation, so we feel that the proof given above might be the best way to proceed.

Under the hypothesis of \ref{Main},  and assuming that $A$ and $B$  are unital,
Pitts  shows that there are
maximal and minimal C*-norms on the *-subalgebra of $B$ linearly spanned by
$\NAB$ \cite[Corollary 7.5]{SRI}.  The completions of this subalgebra under such norms is likely to coincide with the full and reduced groupoid
C*-algebras, respectively, so this would perhaps lead to yet another avenue for proving \ref{Main}.

Recalling that a groupoid C*-algebra is nuclear iff the groupoid is amenable \cite{ClaireJean}, we have the following consequence:

\state Corollary \label Coro
  Under the conditions of \ref{Main}, if $B$ is moreover nuclear, then the groupoid $\G$ mentioned there
is amenable, and hence
  $$
  B\simeq C^*(\G,\Lbpar) \simeq C^*_\mu (\G,\Lbpar) \simeq C^*\red(\G,\Lbpar).
  $$

\Proof
  Since the norm $\mu $ given by \ref{Main} dominates the reduced norm, we get a surjective *-homomorphism
  $$
  B \simeq C^*_\mu (\G,\Lbpar) \to C^*\red(\G,\Lbpar),
  $$
  which implies that $C^*\red(\G,\Lbpar)$ is a quotient of $B$, and hence a nuclear algebra.  This in turn implies that $\G$ is amenable.
\endProof

\endsection

\startsection The opaque ideal

\sectiontitle

The original motivation for pursuing the ideas that led to this paper was the study of the opaque ideal of regular inclusions.  In order to
describe it, let us assume throughout this section that:

\fix $B$ is a C*-algebra, and that
  $$
  A\subseteq B
  $$
  is an abelian,  regular, sub-C*-algebra.  As usual we will denote the spectrum of $A$ by $X$.

Referring to the evaluation states $\varphi _x$ introduced in the proof of \ref{LemmaCondexp}, let $J_x$ be the kernel of
$\varphi _x$, so that $J_x$ is an ideal in $A$ (although quite likely not an ideal in $B$).
  The \emph{opaque ideal} for the inclusion ``$A\subseteq B$'' is defined in \EP{Definition}{11.4}{\DefBlackIdeal} as being
  $$
  \Delta = \bigcap_{x\in X} J_xB.
  $$
  See \EP{Proposition}{11.2.ii}{\IntroBlack.ii} for a proof that the opaque ideal is indeed a closed two-sided ideal in $B$.

It is not
very often that one comes across an example of an inclusion with a nonzero opaque ideal, and in fact in well
behaved cases, such as when $B$ is abelian \cite[Section 4]{EPZ}, or when $B$ is a reduced groupoid C*-algebra
\EP{Proposition}{15.3.i}{\DescribeGray.i},
the opaque ideal is zero.
  For an example of a nonzero $\Delta $ the reader is referred to  \cite[Proposition 4.1]{EPZ}.  The fact that this
example employs an action of the free group causes one to suspect that the lack of amenability of the free group is the main
culprit behind this phenomena.  In fact \cite[Question 4.2]{EPZ} may be thought of as an attempt to rule out bad
behavior of non-amenable groups via the requirement of nuclearity.

\bigskip
  Let us now give a characterization of the opaque ideal in terms of extended pure states.  Recalling that
$\P(A)$ denotes the set of all pure states on $A$, let us define the set of \emph{extended pure states} on $B$ by
  $$
  \Ep = \big \{\psi \in \rsboxtwo{S}\, (B): \psi |_A\in \P(A)\big \},
  $$
  where $\rsboxtwo{S}\, (B)$ is the state space of $B$.

In plain words $\Ep$ consists of all states on $B$ obtained as the extension of some pure state on $A$.  Notice that,
although every pure state on $A$ admits a \emph{pure} extension to $B$, some might also admit a  \emph{mixed} (i.e.~non-pure)
extension. In other words, $\Ep$ might very well  contain mixed states.  An exception is of course the situation in which $A$
has the extension property relative to $B$, in which case $\Ep$ is nothing but the set formed by the $\psi _x$,  all of
which are necessarily pure.

\state  Proposition \label OpaqueIdeal
  \izitem
  \zitem
  If $x\in  X$,  and if $\psi $ is a state on $B$ extending $\varphi _x$, then $\psi $ vanishes on $J_xB$.
  \zitem
  Every $\psi $ in $\Ep$ vanishes on $\Delta $.
  \zitem
  $
  \Delta  = \big \{b\in  B: \psi (b^*b)=0, \text{ for all } \psi \in \Ep\big \}.
  $
  \zitem In case  $A$ has the extension property relative to $B$, we moreover have that
  $$
  \Delta = \big \{b\in B: E(b^*b)=0\big \},
  $$
  where $E$ is the unique conditional expectation given by \ref{LemmaCondexp}.

\Proof
  If $\psi $ extends $\varphi _x$ then,
  given $a\in  J_x$ and $b\in  B$, we have by Cauchy-Schwartz that
  $$
  |\psi (ab)|^2\leq   \psi (aa^*)  \psi (b^*b) =   \langle a^*a, x\rangle \,  \psi (bb^*) = 0,
  $$
  proving  (i),  while (ii) follows immediately from (i).

Regarding (iii), notice that
the inclusion ``$\subseteq $'' follows from (ii).
In order to  prove the reverse inclusion, let $b$ be an element belonging to the set in
the right-hand-side of (iii).   To prove that $b$ lies in $\Delta $, and making use of \EP{Proposition}{11.2}{\IntroBlack},
it suffices to prove that $E_x(b^*b)=0$, for every $x$ in $X$.   Arguing by contradiction, suppose that this is not so for
some $x$.  In
this case one can find a state $\rho $ on the isotropy algebra $B(x)$, such that $\rho (E_x(b^*b))\neq 0$.  By
\EP{Proposition}{3.20}{\CorrespondenceExtendedStates} the composition
  $$
  \psi : B\labelarrow {E_x} B(x) \labelarrow \rho {\bf C}
  $$
  is a state on $B$ vanishing on $J_x$,  and hence  extending  $\varphi _x$,  so that $\psi \in \Ep$.  We then have that
  $$
  0\neq  \rho (E_x(b^*b)) = \psi (b^*b),
  $$
  a contradiction, hence proving (iii).
The last point follows easily from (iii), the fact that
  $$
  \Ep = \{\psi _x:x\in X\},
  $$
  and the definition of $E$ in \ref{LemmaCondexp}.
\endProof

In \cite[Definition 2.4]{SRI},  Pitts introduces the set $\text{Mod}(A, B)$ formed by the  \emph{$A$-modular states}, which,
in our case, coincides with $\Ep$.
  Furthermore, in \cite[Theorem 4.8]{DonsigPitts}, Donsig and Pitts prove that, when $A$ and $B$ are unital, and $A$ has
the extension property relative to $B$, then the subspace of $B$ consisting of the elements annihilated by all modular states
coincides with the left kernel of the conditional expectation $E$, hence the opaque ideal by \ref{OpaqueIdeal.iii}.
They also show that  the quotient of $B$
by said ideal is a C*-diagonal, thus obtaining a characterization of C*-diagonals in terms of the opaque ideal.

We spell out the details of this characterization below, for the sake of completeness, and also to allow for the non-unital case as well.
See \cite[Proposition 8]{Pitts} for similar characterizations.

\state Proposition   \cite[Theorem 4.8]{DonsigPitts}
Let $B$ be a C*-algebra and let $A\subseteq B$ be a regular, abelian sub-C*-algebra.  Then the
following are equivalent:
  \iaitem
  \aitem $A$ is a C*-diagonal in $B$ in the sense of\/ \cite[Definition 1.3]{Kumjian},
  \aitem $A$ satisfies the extension property relative to $B$,  and the opaque ideal vanishes.

\Proof
Assuming (a) we have that $A$ satisfies the extension property relative to $B$ by \cite[Proposition 1.4]{Kumjian}.
Also,  the hypotheses guarantee the existence of a faithful conditional expectation, so the expectation $E$ of
\ref{LemmaCondexp}, being the unique such, must be faithful.  Therefore the opaque ideal vanishes by
\ref{OpaqueIdeal.iv}.

Conversely, by  \ref{Main} we may assume that $B=C^*_\mu (\G,\Lbpar)$, and $A=C_0\big (\Gz\big )$, where
$(\G,\Lbpar)$ is a {\tpeg}.  The composition
  $$
  B = C^*_\mu  (\G,\Lbpar) \labelarrow \pi   C^*\red (\G,\Lbpar) \labelarrow F  C_0\big (\Gz\big ) = A,
  $$
  where $\pi $ is the natural map,  and  $F$ is the standard conditional expectation on the reduced algebra, is clearly a
conditional expectation, hence coincides with $E$ by \ref{LemmaCondexp}.  Since $F$ is well known to be faithful, it
follows that
  $$
  \text{Ker}(\pi ) = \big \{b\in B: F(\pi (b^*b))=0 \big \}  = \big \{b\in B: E(b^*b)=0 \big \} \={OpaqueIdeal.iv} \Delta .
  $$
  The assumption that $\Delta =\{0\}$ then tells us that $B$ coincides with  $C^*\red (\G,\Lbpar)$, so the conclusion follows
from \cite[Theorem 2.9]{Kumjian}.
\endProof

To conclude,  let us give an affirmative answer to \cite[Question 4.2]{EPZ} in the special case that $A$ satisfies  the
extension property.

\state Theorem
  Let $B$ be a C*-algebra and let $A\subseteq B$ be a regular, abelian sub-C*-algebra.  If $B$ is nuclear, and  $A$  has the
extension property relative to $B$,  then $\Delta =\{0\}$.

\Proof
  By \ref{Coro}, we have that
  $B\simeq C^*\red(\G,\Lbpar)$,
  so the conclusion follows from \EP{Proposition}{15.3.i}{\DescribeGray.i}.
\endProof
\endsection

\references

\Bibitem Akemann
  C. Akemann and F. Schultz;
  Perfect C*-algebras;
  Mem. Amer. Math. Soc, 326 (1985)

\Bibitem ClaireJean
  C. Anantharaman-Delaroche and J. Renault;
  Amenable Groupoids;
  L'Enseignement Math\'emati\-que, 2000

\Article JAnderson
  J. Anderson;
  Extensions, restrictions, and representations of states on C*-algebras;
  Trans. Amer. Math. Soc., 249 (1979), 303-329

\Article Archbold
  R. J. Archbold;
  Extensions of Pure States and Projections of Norm One;
  J. Funct. Analysis, 165  (1999), 24-43

\Article ArchboldBunceGregson
  R. J. Archbold, J. W. Bunce and K. D. Gregson;
  Extensions of states of C* -algebras, II;
  Proc. Roy. Soc. Edinb, 92A (1982), 113-122

\Article DonsigPitts
  A. P. Donsig and D. R. Pitts;
  Coordinate systems and bounded isomorphisms;
  J. Operator Theory, 59 (2008), 359-416

\Bibitem ExelPitts
  R. Exel and D. Pitts;
  Characterizing Groupoid C*-algebras of non-Hausdorff \'etale groupoids;
  arXiv:1901.09683 [math.OA], 2019

\Bibitem EPZ
  R. Exel, D. Pitts and V. Zarikian;
  Exotic ideals in represented free transformation groups;
  \hfill\break arXiv:2109.06293 [math.OA], 2021

\Article FeldmanMoore
  J. Feldman and C. Moore;
  Ergodic equivalence relations, cohomology, and von Neumann algebras. I;
  Trans. Amer. Math. Soc., 234  (1977), 289-324

\Article KadisonSinger
  R. V. Kadison and I. M. Singer;
  Extensions of Pure States;
  Amer. J. of Math,  81 (1959),  383-400

\Article Kumjian
  A. Kumjian;
  On C*-diagonals;
  Canad. J. Math., 38 (1986), 969-1008

\Bibitem SRI
  D. R. Pitts;
  Structure for regular inclusions. I;
  J. Operator Theory, 78 (2017), 357-416

\Bibitem Pitts
  D. Pitts;
  Normalizers and approximate units for inclusions of C*-algebras;
  arXiv:2109.00856 [math.OA], 2021.

\Article Renault
  J. Renault;
  Cartan subalgebras in C*-algebras;
  Irish Math. Soc. Bulletin, 61 (2008), 29-63

\Bibitem DarthVader
  {\iss Darth Vader} (user);
  Comment on question {\it `Do extensions of pure states separate points?'\kern2pt};
  \hfill\break {\itt https://mathoverflow.net/questions/405832/do-extensions-of-pure-states-separate-points}

\endgroup
\close